\documentclass{amsart}
\usepackage{amsthm,amssymb}
\newtheorem{theorem}{Theorem}

\newtheorem{lemma}[theorem]{Lemma}

\begin{document}
\title[Polynomials ]{New bounds on the Hermite Polynomials}
    
\author[I. Krasikov]{Ilia Krasikov}
     
\address{   Department of Mathematical Sciences,
            Brunel University,
            Uxbridge UB8 3PH United Kingdom}
\email{mastiik@brunel.ac.uk}
			      
			       
				

\begin{abstract}
We shall establish two-side explicit inequalities,
which are asymptotically sharp up to a constant factor,
on the maximum value of 
$|H_k(x)| e^{-x^2/2}, $
on the real axis,
where $H_k$ are the Hermite polynomials. 
\end{abstract}

\maketitle
\section{Introduction}
We refer to \cite{szego} for the required 
definitions and basic properties of the Hermite polynomials
$H_k(x)$
mentioned in the sequel.
There are a few known upper bounds on the function
$(H_k(x))^2 e^{-x^2},$ see e.g. \cite{abr}, mainly obtained as a specialization
of a more general
case of the Laguerre polynomials.
Some recent results can be found in \cite{mich}.
However, it seems that all presently known inequalities are larger by
the factor $k^{1/6}$
than the true asymptotic value.
The aim of this note is to establish the following two-side bounds, which are sharp up to
constant factors.
\begin{theorem}
\label{mainth}
Let 
$$M_k= (2k)^{1/6} \max_x \left((H_k(x))^2 e^{-x^2} \right),$$
then for $k \ge 6,$
$$\frac{27}{61}  C_k < M_k
< \frac{2}{3} C_k 
\exp{\left( \frac{15}{8} \left(1+ \frac{12}{4(2k)^{1/3}-9}\right) \right)} .
$$
where
$$C_k =\frac{2k \sqrt{4k-2} \; k!^2}{\sqrt{8k^2-8k+3} \;
(k/2)!^2},$$ for $k$ even, and $$C_k=\frac{ \sqrt{16k^2-16k+6} \;
k!(k-1)!}{\sqrt{2k-1} \; ((k-1)/2)!^2}$$ for $k$ odd.
\end{theorem}
The constants 
are not, of course, best possible and can be improved at the cost of
more extensive calculations. 
\\
The asymptotic of the Hermite polynomials in the transition region,
where the maximum of the function $(H_k(x))^2 e^{-x^2}$ is attained,
is given by the following
classical formula \cite{szego}.
For $x=\sqrt{2k+1}-2^{-1/2} 3^{-1/3} k^{-1/6} z,$ z bounded,
$$(H_k(x))^2 e^{-x^2}=3^{2/3} \pi^{-3/2} 2^{k+ \frac{1}{2}} k! \; k^{-1/6}
\left(A(z)+O(k^{-2/3}) \right)^2.
$$
Here $A(z)$ is the Airy function which can be defined
by means of the Bessel functions $J_\nu$
of the first kind,
$$ A(z)=
\frac{ \pi \xi }{3}  \left( J_{- \frac{1}{3}} (2 \xi^3 )+J_{\frac{1}{3}}(2 \xi^3 ) \right),
$$
where $\xi= \sqrt{\frac{z}{3}}.$
\\
The absolute maximum of $A(z)$ is attained for $z=1.46935...,$ and is equal to $1.1668...$.
Therefore, asymptotically, the optimal $x$ is $\sqrt{2k}-1.61723 (2k)^{-1/6},$
and
$$\lim_{k \rightarrow \infty } \frac{M_k}{C_k} =0.715452 ...$$
Thus, for
sufficiently large $k,$ our upper bound ($\frac{2}{3} e^{15/8}=4.347...$)
is about six times greater
than the corresponding asymptotic value.
\\
We want to stress that in principal
(yet, it would require a formidable amount of calculations) our method combining
the ideas of \cite{fk,kdif,kz} can be applied to
a more general case of the Laguerre and Jacobi polynomials
with parameters arbitrarily depending on $k.$ 
In such a general situation we don't even know the corresponding
asymptotics, see e.g. \cite{erdel}.
\section{Proofs}
To prove Theorem \ref{mainth} we will need to
perform a substantial amount of straightforward but rather tedious calculations.
We used Mathematica to handle them.
\\
As $(H_k(x))^2$ is an even function
it is enough to consider the case $x \ge 0,$
and this will be assumed in the sequel.
For a fixed $k$ 
it will be convenient to introduce two functions $y=y(x)=2k-x^2,$
and the 
logarithmic derivative $t=t(x)=H'_k(x)/H_k(x).$
\\
We also make use 
of the differential equation
\begin{equation}
\label{difequ}
f''=2x f'-2k f, \, \, \, f=H_k(x),
\end{equation}
in a pure algebraic manner to exclude all the derivatives
of $H_k$ of order greater than one when such appear.
\\
Let $x_{1k}<...<x_{kk}$ be the zeros of $H_k(x).$
It is well known that
that the successive relative maxima of $H_k(x) e^{-x^2/2}$
form an increasing sequence for $x \ge 0,$ \cite{szego}. 
Observe that $\frac{d}{dx} H_k(x) e^{-x^2/2} =0,$ just means
$t(x)=x.$  A quick inspection of the graph of the function $t(x),$
consisting of $k+1$ decreasing branches, and 
the corresponding straight line reveals that their last intersection
occurs for  some $x = \omega >x_{kk}.$
Thus, the absolute maximum of the function
$ (H_k(x))^2 e^{-x^2}$ is attained for $x = \omega .$
\\
We will deduce Theorem \ref{mainth} from the following bounds established in \cite{fk}.
\begin{theorem}
\label{th2}
For $x^2 < 2k-\frac{3}{2},$  
$$
(H_k(x))^2 e^{-x^2} \le C_k
F_k(y) G_k (y),
$$
where
$$
F_k(y)=\frac{2y^2-4y+3}{\sqrt{y(4y^4-12y^3+9y^2+10k y-12k)}},
$$
$$
G_k(y)= \exp{(\frac{15(2k-y)}{2y(2y-3)^2}) },
$$
Moreover the inequality is sharp  in a sense that
$$
(H_k(x)^2 e^{-x^2} \ge C_k
F_k(y)/ G_k (y),
$$
for all
the roots of the equation
\begin{equation}
\label{interl}
x y(2y-3)H_{k}(x)=(2y^2-4y+3)H_{k-1} (x),
\end{equation}
that is
at a point between any two
consecutive zeros of $H_k .$
\end{theorem}
We need the following technical result.
\begin{lemma}
\label{incr}
For $k \ge 2,$
the function $v_1(y)=F_k(y)G_k(y)$ decreases in $y $ 
for $ y \ge 2.$
The function $v_2(y)=F_k(y)/G_k(y)$ decreases in $y$ for $y \ge 3 \; (2k)^{1/3}.$
\end{lemma}
\begin{proof}
It is easy to check by the substitution $y:=y+2,$ that 
for  $y \ge 2,$
the denominator of $F_k(y)$
is positive.
Now, $v_1$ is a decreasing function since for $k \ge 2, y \ge 2,$
$$G_k^{-2}(y) \frac{dv_1}{dy} =
- \; \frac{(2y^2-4y+3) A(k,y)}{y^3(2y-3)^3(4y^4-12y^3+9y^2+10k y-12k)^2} <0,$$
where $A(k,y)$ is a polynomial with $deg_k A=2, \, deg_y A=10.$
Indeed one can check that $A(k,y) >0,$
as it is transformed into a polynomial with only
positive terms by the substitutions $k:=k+2, \, y:=y+2.$
\\
Similarly, for $v_2$ we have
$$G_k^2(y) \frac{d v_2}{dy}=
- \; \frac{(2y^2-4y+3) B(k,y)}{y^3(2y-3)^3(4y^4-12y^3+9y^2+10k y-12k)^2} <0,$$
where $B(k,y)$ is a polynomial with $deg_k B=2, \, deg_y B=10.$
It is left to check that
the substitution 
$$y:=y+3 \; (2k)^{1/3}, \, \, k:=k+2,$$  
transforms $B(k,y)$ it into a polynomial with only
positive terms.
We omit the details.
\end{proof}
Let $g=g(x)$ be a real polynomial with only real zeros $x_1,...,x_k.$
The following inequality is called the Laguerre inequality,
\begin{equation}
\label{lagineqv}
\frac{{g'}^2-g g''}{g^2}= \sum_{i=0}^k \frac{1}{(x- x_i)^2} >0,
\end{equation}
and will be our main technical tool in this note.
It is worth noticing that Theorem \ref{th2} was established by
applying a higher order generalization of (\ref{lagineqv}).

To simplify the formulas in the sequel
we will use the substitution $k=\frac{27 m^{12}-1}{54m^6},$ that is
$m=\left( k+\sqrt{k^2+ \frac{1}{27}} \, \, \right)^{1/6}.$
\begin{lemma}
\label{xt}
$\omega < (m^4-\frac{1}{3})^{3/2} m^{-3}.$
\end{lemma}
\begin{proof}
We put $f(x)=H_k(x),$ $g(x)=f(x)-f'(x)/q,$ where $q$ is a parameter independent on $x.$
Observe that $g(x)$ has only real zeros as well.
Therefore, by (\ref{lagineqv}),
$$
U(t,x,q)=\frac{q^2 ({g'}^2-g g'') }{f^2} \ge 0,
$$
Using
(\ref{difequ})
we get
$$
U(t,x,q)=(2k+q^2 -2-2q x)t^2-2(2k x+q^2 x-2q x^2-q)t+2k(2k+q^2-2q x) \ge 0,
$$
and for $t=x,$
\begin{equation}
\label{cub}
U(x,x,q)=2q x^3-(2k+2+q^2)x^2-2q(2k-1)x+2k(2k+q^2) \ge 0.
\end{equation}
First we shall show that this inequality implies $x <\sqrt{2k-1}.$
Choosing $q=0,$ we have
$U(x,x,0)=2(2k^2-(k+1)x^2) \ge 0,$
hence 
$$x \le \sqrt{\frac{2k^2}{k+1}}<\sqrt{2k-1}.$$
The optimal value of 
$q=q_0$
corresponds to the case when the discriminant $\Delta$ of
(\ref{cub}), considered as a quadratic in $q,$ is zero. That is, practically, one has to solve the system
$U(x,x,q_0)=0, \, \Delta=0,$
yielding 
$$x=(m^4-\frac{1}{3})^{3/2} m^{-3}, \, \, q=q_0= \frac{(3m^4-3m^2-1) \sqrt{9m^4-3}}{9m^3}.$$
Formally this can be justified as follows.
The cubic equation $u(x)=U(x,x,q_0)=0,$ has three real zeros 
$$x_1 <0 <x_2=(m^4-\frac{1}{3})^{3/2} m^{-3} < \sqrt{2k-1} <x_3.$$
Indeed $u(-\infty )=-\infty ,\, \,u(\infty )=\infty ,$ and
$$
u(0)=2k(2k+ q_0^2) >0, \, \, \, u( \sqrt{2k-1} \, )=
- \, \frac{2(9m^8-15m^4+1)}{9m^4} <0 .
$$
Thus we conclude that $t=x,$ implies 
$$x<x_2 =(m^4-\frac{1}{3})^{3/2} m^{-3}.$$
\end{proof}
Now we can prove the upper bound of Theorem \ref{mainth}.
\begin{lemma}
\label{voz}
Let $k  \ge 6 ,$ then
$$
M_k <
\frac{2}{3} C_k (2k)
\exp{\left( \frac{15}{8} \left(1+ \frac{45}{4(2k)^{1/3}-9}\right) \right)}
$$
\end{lemma}
\begin{proof}
Observe that $k \ge 6,$ implies $m >3/2.$
By the previous lemma 
$$y(\omega ) =2k- \omega^2 > m^2-\frac{1}{3m^2} =y_0 >2.$$
Thus, the conditions of Lemma \ref{incr} are fulfilled and it is left to
estimate $F_k ( y_0 )G_k (y_0).$
The required upper bound on $G_k(y_0)$ follows from 
$$
\frac{2k-y_0}{ y_0 (2 y_0-3)^2 }= 
\frac{(3m^4-1)^2}{(6m^4-9m^2-2)^2} \le \frac{1}{4}\left(1+ \frac{12}{4m^2-9} \right) 
<\frac{1}{4}\left(1+ \frac{12}{4(2k)^{1/3}-9} \right).
$$
\\
Straightforward calculations also yield
$$F_k(y_0)<\frac{2}{3 m} <\frac{2}{3} (2k)^{-1/6},$$ 
and the result follows.
\end{proof}
To demonstrate the lower bound of Theorem \ref{mainth} we observe that by $H'_k(x)=2k H_{k-1}(x),$
equation (\ref{interl})
can be rewritten as
$$
t= \frac{2k x y(2y-3)}{ 2y^2-4y+3}.
$$
The right hand side is a continuous odd function, positive for $x>0,$ hence intersecting
all the branches of $t.$
Unfortunately, the intersections points to the right of $x_{kk}$
violate the restrictions of Lemma \ref{incr}.
Therefore we will choose an intersection point $\tau$ from the interval
$(x_{k-1,k},x_{kk}).$ In fact $\tau$ is greater than the largest zero of
$H'_k(x)$, that is $x_{k-1,k-1},$ as $t(\tau )>0.$
Thus, for the lower bound we just calculate $F_k(x_{k-1,k-1})/G_k (x_{k-1,k-1})$
and show that $2k-x_{k-1,k-1}^2 >3 \; (2k)^{1/3},$ 
the condition imposed by Lemma \ref{incr}.
The last claim is justified by the following lemma which maybe of independent interest. 
\begin{theorem}
\label{ozkor}
For $k >2$ the largest zero $x_{kk}$ of $H_k(x)$ satisfies
$$ x_{kk} > \sqrt{2k}-\frac{9}{4} (2k)^{-1/6} $$
\\
In particular,
$2k-x_{k-1,k-1}^2 >3 \; (2k)^{1/3},$ for $ k \ge 3.$
\end{theorem}
\begin{proof}
We use the method of \cite{kz} based on the so-called Bethe ansatz equations.
First we shall prove that for any $x > x_{kk},$
\begin{equation}
\label{eqd}
2k-x^2 < \frac{1}{(x-x_{kk})^2} + \frac{2k-2-x_{kk}^2}{3}.
\end{equation}
Using the differential equation (\ref{difequ})
to exclude higher derivatives we get
$$\frac{{f'}^2-f f''}{f^2} =t^2-2x t+2k=(t-x)^2+2k-x^2 \ge 2k-x^2.$$
On the other hand, by (\ref{lagineqv}) 
$$\frac{{f'}^2-f f''}{f^2}=\sum_{i=1}^k \frac{1}{(x-x_{ik})^2} .$$
Therefore, for any $x > x_{kk},$
$$2k-x^2 \le \sum_{i=1}^k \frac{1}{(x-x_{ik})^2} < \frac{1}{(x-x_{kk})^2} +
\sum_{i=1}^{k-1} \frac{1}{(x_{kk}-x_{ik})^2}.$$
The last sum is equal to
$$
\frac{2k-x_{kk}^2-2}{3},
$$
and can be easily calculated as the limit
$$
\lim_{x \rightarrow x_{kk}} \left(\frac{{f'}^2-2x f' f+2k f^2}{f^2} -\frac{1}{(x-x_{kk})^2}
\right),
$$
by applying 
L'H\^{o}pital's rule
four times
and substituting $f''$ from (\ref{difequ}) at each step.
This yields (\ref{eqd}).
\\
Put now $k=s^6/2, \,
x=x_{kk}+1/s. $
Then (\ref{eqd}) can be rewritten as
$$
2s^2 x_{kk}^2+6s x_{kk}-2s^8+3s^4-2s^2+3 >0.
$$
Hence
$$x_{kk} > \frac{\sqrt{4s^8-6s^4+4s^2+3}-3}{2s} > s^3-\frac{9}{4s},$$
and the result follows.
\\
The second claim is a matter of simple calculations.
\end{proof}
It is worth noticing that the obtained result is quite precise, as 
$$x_{kk}< \sqrt{2k+1}-6^{-1/3} (2k+1)^{-1/6}i_1 , $$
where $i_1$ is the least positive zero of the Airy's function,
$(6^{-1/3} i_1=1.85574...)$ \cite{szego}.
\\
Now we are in the position to complete the proof of Theorem \ref{mainth}.
\begin{lemma}
\label{max}
$$M_k > \frac{27}{61} C_k .$$
\end{lemma}
\begin{proof}
We just have to calculate the value of the function
$ \frac{ (2k)^{1/6} F_k (y)}{G_k(y)}$
for
$$x= \sqrt{2k-2}-\frac{9 }{4} (2k-2)^{-1/6}.$$
One can check that it
has the only minimum $0.44265... > 27/61,$ for $k=46,$
(notice that the asymptotic value, for $k \rightarrow \infty ,$ is $0.4586...,$
and only
slightly better).
This completes the proof.
\end{proof}


\begin{thebibliography}{8}
  
\bibitem{abr}
M.\,Abramowitz, I.A.\,Stegun, Handbook of Mathematical
Functions, Dover, New York, 1964.
\bibitem{erdel} T.\,Erdelyi, A.P.\,Magnus, P.\,Nevai,
{\em Generalized Jacobi weights, Christoffel functions, and
Jacobi polynomials}, SIAM J. Math. Anal. 25 (1994), 602-614.
\bibitem{fk}
W.H. Foster and I. Krasikov,
Explicit bounds for Hermite  polynomials in the oscillatory region,
{\em LMS J. Comput. Math.}, Vol.3 (2000) 307-314.
\bibitem{kdif}
I. Krasikov, On zeros of polynomials and allied functions satisfying
second order differential equation,
{\em East J. Approx.}, 9 (2003) 51-65.
\bibitem{kz}
I. Krasikov, On extreme zeros of classical orthogonal polynomials, submitted.
\bibitem{mich}
M.\,Michalska, J.\,Szynal, A new bound for the Laguerre polynomials,
{\em J. Comput. Appl. Math.} 133 (2001) 489-493.
\bibitem{szego} 
G.\,Szeg\"{o}, {\em Orthogonal Polynomials},
Amer. Math. Soc. Colloq. Publ., v.23, Providence, RI, 1975.
     
\end{thebibliography}
\end{document}